\documentclass[a4paper,11pt]{article}

\usepackage{amssymb, amsmath, eucal}

\newtheorem{teorema}{Theorem}

\newtheorem{lema}[teorema]{Lemma}
\newtheorem{proposicao}[teorema]{Proposition}
\newtheorem{corolario}[teorema]{Corollary}

\newenvironment{prova}{\setlength{\parindent}{0pt}\textbf{Proof.}}{\hspace{\stretch{1}}$\Box$}
\newenvironment{provap}{\setlength{\parindent}{0pt}\textbf{Proof of proposition.}}{\hspace{\stretch{1}}$\Box$}

\title{Paraconsistentization and many-valued logics}

\author{Edelcio G. de Souza$^{1}$ \\
Alexandre Costa-Leite$^{2}$ \\ 
Diogo H. B. Dias$^{3}$\\
\footnotesize $^{1}$ University of Sao Paulo (BR) \ $^{2}$ University of Brasilia (BR) \\ \footnotesize $^{3}$ State University of Northern Parana (BR)}

\date{}

\begin{document}

\maketitle

\begin{abstract}
This paper shows how to transform explosive many-valued systems into
paraconsistent logics. We investigate mainly the case of three-valued
systems exhibiting how non-explosive three-valued logics can be obtained
from them. 
\end{abstract}

\section{Introduction}

An explosive logic can be transformed into a non-explosive one by
means of methods of \emph{paraconsistentization}. There are many
ways which can be used to perform this task of converting explosive
into non-explosive logics. In general, it is
theoretically possible to \emph{paraconsistentize} all non-paraconsistent
systems.

In the history of paraconsistency, we can find, for instance, three accounts
which can be used to introduce a paraconsistent dimension into an explosive system.  These methods were, indeed, classified as \emph{concrete} forms of paraconsistentization in \cite{egs-acl-dd-3}.
The first one developed by Stanis\l aw  Ja\'skowski in \cite{jaskowski} focuses in a way to define 
a paraconsistent discussive logic from a standard modal logic, while the second one proposed 
by Newton da Costa in \cite{dacosta} introduces
a paraconsistent negation using the notion of \emph{well-behavior} in such a way that new 
properties of a weaker negation not satisfying \emph{ex falso} are generated (cf. also \cite{dacosta11}) . The third approach suggested by Graham Priest
in \cite{priest0} considers a new logical designated value besides truth
and this gives rise to a paraconsistent logic. These three 
accounts can be seen as methods for paraconsistentizing classical logic and, 
in some sense, they are mechanisms of paraconsistentization which can be applied to a wide range of logics.

All previous techniques are not unified by a standard strategy to transform
a given logic into a paraconsistent one. Indeed, each one uses a particular procedure to formulate paraconsistency in an specific level. The initial concept and idea of \emph{paraconsistentization} 
has been proposed in \cite{costaleite}.  Afterwards, in two recent papers (cf. \cite{egs-acl-dd-1} 
and \cite{egs-acl-dd-2}), it has been 
shown how to turn a given logic into a paraconsistent one by a very precise locally developed methodology which works for a great variety of systems. In \cite{egs-acl-dd-1},
the idea of \emph{paraconsistentization} is presented in the context of category theory and
it uses abstract logic as a main source.  In this way, paraconsistentization appears as an endofunctor in the category of logics preserving some basic properties of the initial logic. In \cite{egs-acl-dd-2}, considering notions
of \emph{axiomatic formal systems} and \emph{valuation structures}, the quest of paraconsistentization
is introduced by means of the concept of \emph{paradeduction} in axiomatic formal
systems and \emph{paraconsequence} in valuation structures. Thus, it is possible to
paraconsistentize proof systems and semantics showing how some properties
are invariant under paraconsistentization. In this particular case, this method
preserves soundness and completeness.

The above papers settled the basic theory of paraconsistentization, up to this level.
Despite the fact the there are many unknown methods of paraconsistentization, this article shows how to turn some explosive many-valued logics into paraconsistent ones  using the basic idea initially proposed in  \cite{egs-acl-dd-1}, that is:
given a logic $L= \langle For, \vDash_{L} \rangle$, the paraconsistentization of 
$L$ is a logic given by $\mathbb{P}(L) = \langle For, \vDash^{\mathbb{P}}_{L} \rangle$ such that: $\Gamma \vDash^{\mathbb{P}}_{L} \alpha$ if, and only if, 
$\text{there exists }  \Gamma' \subseteq \Gamma, L\text{-consistent such that  } \Gamma' \vDash_{L} \alpha$.  This has been classified, besides other ideas as those of \cite{bensusan-greg}, as a type of \emph{abstract} paraconsistentization.

The technique of using consistent sets departing from some forms of contradictory sets to manage inconsistency has a long history  
which starts at least since the works of Nicholas Rescher and Ruth Manor in \cite{rescher-manor}. In this approach, consequence relations are generated by means of the notion of \emph{maximal consistent sets}. We
have addressed what are the main differences between this and our approach in \cite{egs-acl-dd-2}.
Afterwards, studies and researches somewhat related with this initial departure point can be mentioned: \cite{amgoud, arieli2, avron, benferhat, brown, grant-sub, subrahmanian, vandeputte}. In
the spirit of preservationism, which studies what properties are preserved by a consequence
relation generated by a kind of forcing, there is also a type of general paraconsistentization developed by Gillman Payette in \cite{payette}. In spite of using consistent sets to generate a paraconsistent version of the initial consequence relation, our main target is not
to explore how reasoning with consistent sets works. Instead, we are, indeed, interested in the
logical movement, that is, we are mainly motivated by converting a logic not able to deal with 
inconsistencies into a logic able to tolerate, control and manipulate inconsistent inferences (in a
close relation with approaches such as those of \cite{dacosta-vernengo} and \cite{desouza5}).
Our idea is to study the process of generating paraconsistent versions of previous logics, and
this can be done in many different directions, not only using consistent sets.

In this paper, we deal mainly with systems $L_{3}$, $G_{3}$, $K_{3}$ and $LP$ defined by means of logical matrices, though it is obviously
possible to extend the same approach to the whole hierarchies $L_{n}$ and $G_{n}$. We show what are the properties preserved or lost by the paraconsistentized version of these systems. This paper attempts to explore the
universe of paraconsistent many-valued logics whose relevance appears when we take
into account motivations and basic intuitions of each of the many-valued systems considered.\footnote{For essential information on motivations, relevance, applications and aims of many-valued logics, there
are some substantial available literature, see \cite{gottwald}, \cite{malinowski}, \cite{epstein}, \cite{bolc} and \cite{priest1}. Concerning three-valued paraconsistent logics it is worth mentioning \cite{arieli-avron}, \cite{arieli-avron-zamanski} and \cite{coniglio1}.} Furthermore, for
us, these paraconsistentized many-valued systems can play a substantial improvement in our notion of  \emph{partial justifications} (cf. \cite{dacosta1}), as this notion essentially allows the coexistence of contradictory situations. For instance, a sentence and its negation,
especially those empirical sentences about reality, can be both partially justified. What we have in mind is the following situation. Let $p$ be an atomic sentence. Sometimes, in 
the realm of empirical sciences, it is possible to find justifications for a given sentence
$p$ and also for its negation. These are justifications called by da Costa in \cite{dacosta1} as weak justifications (i.e. partial justifications). He shows that depending of the properties the notion of justification
has, a contradictory situation is generated. We conjecture that using paraconsistentized 
many-valued logics helps in the way we deal with these justifications. So,
interpreting partial justification in the realm of paraconsistentized many-valued systems
can be useful to manipulate scenarios in which we need to handle a justification of a sentence
and a justification of its negation (cf. \cite{costaleite3} and \cite{costaleite4} for a preliminary approach to the problem of partial justifications). 

First, in the logic $L_3$, conceived by \L ukasiewicz in \cite{lukas}, it is feasible to define modal operators such as \emph{possibility}, \emph{necessity} and even the notion of \emph{indetermination} (cf. \cite{malinowski}). These can be used to model forms of weak justification and so, for this reason, paraconsistentizing $L_3$ increases our ability to deal with inconsistent justifications.
Second, $K_3$ was proposed to deal with types of partial information with respect to the truth value of propositions in such a way that Kleene in \cite{kleene} suggests that the intermediate truth-value can be evaluated as ``undefined'' or as ``unknown''. So, this third logical value can be seen as absence of justification and, for this reason, it is also connected with partial justifications. Given this hermeneutics,  a proposition with an undefined truth-value might became true or false (depending of the amount of justification), according to future information which will be added (i.e according to some new justification which can be included to support the truth of a given statement).  But since explosion holds in these systems, contradictory information or justifications are \emph{a priori} excluded. By paraconsistentizing these logics, we provide a formal framework in which these scenarios - contradictory justifications - can be investigated. In this sense, paraconsistentization of logics widens the scope of investigation of these logics. Third, G\"odel intermediate logics have been
developed with the aim of showing that intuitionistic logic cannot be defined
using ``finitely many elements (truth values)'' (G\"odel, page 225, \cite{godel}). G\"odel discovered an infinite hierachy of consistent systems between
intuitionistic and classical logic. Paraconsistentizing this hierarchy shows
that there exists a whole new hierarchy of paraconsistent many-valued systems which does not need, necessarily, to be consistent. The problem of determining where is this new hierarchy (between which logics?) is still open. Moreover, G\"odel's three valued system can also
be applied to the problem of partial justifications. 

It is important to stress two methodological dimension of this proposal. First, since paraconsistentization preserves consistent sets of formulas from the original logic, if no contradiction arises, then the paraconsistent version of a logic is equivalent to the original one. Second, paraconsistentization merely allows reasoning with inconsistent sets, but does not impose a specific interpretation of this inconsistency. In particular, paraconsistentization is neutral with respect to some philosophical problems like those concerning the existence (or not) of true (or real) contradictions. 

In what follows, we apply the standard methodology of \emph{paraconsistentization by consistent sets}
to particular many-valued systems. Using this technique we are able to build paraconsistent
many-valued systems on a very large scale. 

\section{Preliminaries}

These preliminaries establish main terminology and concepts that are used
throughout the text. Standard notions and results on many-valued logics can
be found in \cite{gottwald}, \cite{malinowski}, \cite{epstein}, \cite{bolc} and \cite{priest1}. 

Let us consider an usual \emph{propositional language} $L$ with $\neg$ (negation symbol), $\vee$ (disjunction symbol), $\wedge$ (conjunction symbol), $\rightarrow$ (implication symbol) and propositional letters: $p, q, r,...,p_{1}, q_{1}, r_{1},...$ and so on. The set of propositional letters is denoted by $Prop$ and the set of formulas of $L$, defined as usual, is
denoted by $For$. We use $\alpha$, $\beta$, $\gamma$,..., and $\Gamma$, $\Delta$,... as syntactical variables for formulas and sets of formulas, respectively.

A \emph{(logical) matrix} for $L$ is a 6-tuple
\begin{center}
$M=\langle Val, D, f_{\neg}, f_{\vee}, f_{\wedge}, f_{\rightarrow} \rangle$
\end{center}
such that $Val$ is a non-empty set, $D$ is a non-empty proper subset of $Val$, $f_{\neg}$ is an unary function $f_{\neg}: Val \rightarrow Val$  and $f_{\vee}, f_{\wedge}, f_{\rightarrow}$
are binary functions of the type $f_{\vee}, f_{\wedge}, f_{\rightarrow}: Val \times Val \rightarrow Val$. Elements of $Val$ are called \emph{truth-values} (or simply \emph{values}) and elements of $D$ are called \emph{designated values}.

A $M$-\emph{valuation} for  $L$ is a function 
$$
v: Prop \rightarrow Val.
$$
Every $M$-valuation
 $v$ can be extended, in an unique way, to all elements of $For$ by the following recursive clauses:

\medskip
i. $v(\neg \alpha)= f_{\neg}(v(\alpha))$;

ii. $v(\alpha \vee \beta)= f_{\vee}(v(\alpha), v(\beta))$;

iii. $v(\alpha \wedge \beta)= f_{\wedge}(v(\alpha), v(\beta))$;

iv. $v(\alpha \rightarrow \beta)= f_{\rightarrow}(v(\alpha), v(\beta))$.

\medskip
We denote by $\mathcal{V}_{M}$ the set of all $M$-valuation
for $L$. Let $\Gamma$ be a subset of $For$ and $M$ a matrix for $L$. An element $v \in \mathcal{V}_{M}$ is a $M$-\emph{model} of $\Gamma$ iff (if and only if) $v(\gamma) \in D$, for all $\gamma \in \Gamma$. We denote by $Mod_{M}(\Gamma)$ the set of all $M$-models for $\Gamma$. If $\Gamma = \{\alpha\}$ is an unitary set, we denote the set $Mod_{M}(\{\alpha\})$ by $Mod_{M}(\alpha)$. Let $\Gamma$ be a subset of $For$ and $\alpha$ an element of $For$. We say that $\alpha$ is a $M$-\emph{consequence} of $\Gamma$, in symbols $\Gamma \vDash_{M} \varphi$, if every $M$-model of $\Gamma$ is a $M$-model of $\alpha$, that is, $Mod_{M}(\Gamma) \subseteq Mod_{M}(\alpha)$. We have the following immediate properties:

\medskip
I. If $\alpha \in \Gamma$, then $\Gamma \vDash_{M} \alpha$;

II. If $\Gamma \vDash_{M} \alpha$, then $\Gamma \cup \Delta \vDash_{M} \alpha$;

III. If $\Gamma \vDash_{M} \alpha$ and $\Delta \vDash_{M} \gamma$, for all $\gamma \in \Gamma$, then  $\Delta \vDash_{M} \alpha$.

\medskip
If $\Gamma$ is a set of formulas, the set of $M$-\emph{consequences} of $\Gamma$, denoted by $Cn_{M}(\Gamma)$, is such that:
$$
Cn_{M}(\Gamma) := \{\alpha \in FOR: \Gamma \vDash_{M} \alpha\}.
$$
$Cn_{M}$ can be seen, therefore, as an operator on $\wp (For)$, the set of all subsets of $For$, and $Cn_{M}$ satisfies the Tarskian axioms (cf. \cite{tarski} and \cite{tarski1}) :

\medskip
I'. $\Gamma \subseteq Cn_{M}(\Gamma)$;

II'. $Cn_{M}(\Gamma) \subseteq Cn_{M}(\Gamma \cup \Delta)$;

III'. $Cn_{M}(Cn_{M}(\Gamma))= Cn_{M}(\Gamma)$.

\medskip
 If $\alpha$ is an element of $For$, we say that $\alpha$ is a $M$-\emph{tautology} iff $Mod_{M}(\alpha) = \mathcal{V}_{M}$, that is, every $M$-valuation $v$ is such that $v(\alpha) \in D$. Thus, every $M$-valuation is a $M$-model of $\alpha$. We use the notation $\vDash_{M} \alpha$ to indicate that $\alpha$ is a $M$-tautology. It is easy to see that: $\vDash_{M} \alpha$ iff $\emptyset \vDash_{M} \alpha$. Therefore, the set of all $M$-tautologies is the set $Cn_{M}(\emptyset)$. By monotonicity, the property (II') above, for all $\Gamma \subseteq For$, $Cn_{M}(\emptyset)\subseteq Cn_{M}(\Gamma)$, that is, the set of $M$-consequences of a set $\Gamma \subseteq For$ always contains all the $M$-tautologies.

\section{Paraconsistentization of logics}
We start by reviewing some notions in the domain of paraconsistentization and we follow the presentations in \cite{egs-acl-dd-1} and \cite{egs-acl-dd-2}. 

A pair $L = \langle X,Cn_{L}\rangle$ is a \emph{consequence structure} iff $X$ is a non empty set and $Cn_{L}$ is a mapping in $\wp (X)$. Note that, in general, we use a consequence relation with
no axioms. In a consequence structure of the type $L=\langle X, Cn_{L} \rangle$, $\Gamma \subseteq X$ is $L$-\emph{consistent} if and only if $Cn_{L}(\Gamma) \neq X$. Otherwise, $\Gamma$ is said to be $L$-\emph{inconsistent}. We say that a consequence structure $L = \langle X,Cn_{L}\rangle$ is \emph{normal} iff $Cn_{L}$ satisfies the properties (I')-(III') above. Moreover, we say that $L$ is \emph{compact} iff every $L$-inconsistent set has a $L$-inconsistent finite subset. If $L=\langle X, Cn \rangle$ is a consequence structure, we define a  $\mathbb{P}$-transformation of $L$, called a \emph{paraconsistentization} of $L$, as a consequence structure
$\mathbb{P} (L) = \langle X, Cn^{\mathbb{P}}_{L}\rangle$ such that\footnote{Notice that the domains of $L$ and  $\mathbb{P}(L)$ are the same set $X$.}: For all subset $\Gamma$ of $X$ we have
$$
Cn^{\mathbb{P}}_{L}(\Gamma) = \bigcup \{ Cn_{L}(\Gamma') \in \wp (X): \Gamma' \subseteq \Gamma \text{ and } \Gamma' \text{ is } L \text{-consistent}\}.
$$
Therefore, we have that $\alpha \in Cn^{\mathbb{P}}_{L}(\Gamma)$ iff there exists  $\Gamma' \subseteq \Gamma$, $L$-consistent, such that $\alpha \in Cn(\Gamma')$.

Now, consider a matrix $M$ for $L$. Thus, $M$ yields a consequence structure $\langle For, Cn_{M}\rangle$ as seen above. (We also will use $M$ for the pair $\langle For, Cn_{M}\rangle$, in order to fit the notation above.) In this way, the paraconsistentization ($\mathbb{P}$-transformation) of $M = \langle For, Cn_{M}\rangle$ is a consequence structure $\mathbb{P}(M) = \langle For, Cn^{\mathbb{P}}_{M}\rangle$ such that: For all $\Gamma \subseteq For$ and $\alpha \in For$, we have that:
\begin{center}
$\Gamma \vDash^{\mathbb{P}}_{M} \alpha$ iff there exists $\Gamma' \subseteq \Gamma$, $M$-consistent, such that $\Gamma' \vDash_{M} \alpha$.
\end{center}
Note that if $M$ is a matrix for $L$ such that $M= \langle For, Cn_{M}\rangle$ is the consequence structure associated to $M$, then we have: If $\Gamma \subseteq For$ is $M$-consistent, then $Mod_{M}(\Gamma) \neq \emptyset$. Because, if $Mod_{M}(\Gamma) = \emptyset$, then $Mod_{M}(\Gamma) \subseteq Mod_{M}(\alpha)$, for all $\alpha \in For$. Therefore, $\Gamma \vDash_{M} \alpha$, for all $\alpha \in For$, that is, $Cn_{M}(\Gamma) = For$. Thus, $\Gamma$ is $M$-inconsistent. Otherwise, the converse is valid only in special cases.

\begin{lema}
\label{consistencia}
Let $M=\langle Val, D, f_{\neg}, f_{\vee}, f_{\wedge}, f_{\rightarrow} \rangle$ be a matrix such that $f_{\neg}$ satisfies the following condition:
$$
(*) \ \ \ \ \ \ \text{if } x \in D, \text{ then } f_{\neg}(x) \notin D.
$$
Then, we have that if $Mod_{M}(\Gamma) \neq \emptyset$, then $\Gamma$ is $M$-consistent.
\end{lema}

\begin{prova}
Suppose that $Mod_{M}(\Gamma) \neq \emptyset$ and let $v \in Mod_{M}(\Gamma)$. Suppose, on the contrary, that $\Gamma$ is $M$-inconsistent. Then, $Cn_{M}(\Gamma) = For$, that is, $Mod_{M}(\Gamma) \subseteq Mod_{M}(\alpha)$, for all $\alpha \in For$. Thus,
$$
Mod_{M}(\Gamma) \subseteq \bigcap_{\alpha \in For} Mod_{M}(\alpha).
$$
Let $p$ be propositional letter. Since $f_{\neg}$ satisfies (*), we have $Mod_{M}(p) \cap Mod_{M}(\neg p) = \emptyset$. Therefore,
$$
Mod_{M}(\Gamma) \subseteq \bigcap_{\alpha \in For} Mod_{M}(\alpha) = \emptyset.
$$
But, $v \in Mod_{M}(\Gamma)$ (contradiction!).
\end{prova}

\medskip
In \cite{egs-acl-dd-1}, we have presented a sufficient condition to the $\mathbb{P}$-transformation of a given logic $L$ to be a paraconsistent logic. We have to consider the following essential concepts (see \cite{egs-acl-dd-1}, definition 4.1, p.246). So, let $L = \langle X, Cn_{L}\rangle$ be a consequence structure. We assume that $X$ is endowed with a negation operator $\neg$. Then, we say that:

\medskip
(a) $L$ is \emph{explosive} iff for all $A \subseteq X$, if $x \in X$ is such that $\{x, \neg x\} \subseteq Cn_{L}(A)$, then $A$ is $L$-\emph{inconsistent} ($Cn_{L}(A)=X$). Otherwise, $L$ is called \emph{paraconsistent};

(b) $L$ satisfies \emph{joint consistency} iff there exists $x \in X$ such that $\{x\}$ and  $\{\neg x\}$ are both $L$-consistent but $\{x, \neg x\}$ is $L$-inconsistent;

(c) $L$ satisfies \emph{conjunctive property} iff for all $x, y \in X$, there is a $z \in X$ such that $Cn_{L}(\{x, y\})=Cn_{L}(\{z\})$.

\begin{teorema}
If a consequence structure  $L=\langle X, Cn_{L}\rangle$ is normal, explosive, satisfies joint consistency and also the conjunctive property, then $\mathbb{P}(L)$ is paraconsistent.
\end{teorema}

\begin{prova}
See \cite{egs-acl-dd-1}, theorem 4.2, p.246.
\end{prova}

\begin{proposicao}
\label{prop36}
If a consequence structure  $L=\langle X, Cn_{L}\rangle$ is normal and there is $\alpha \in For$ such that $\{ \alpha \}$ is $L$-inconsistent, then in $\mathbb{P}(L)$ there is no $\mathbb{P}(L)$-inconsistent sets.
\end{proposicao}

\begin{prova}
See \cite{egs-acl-dd-1}, proposition 3.6, p.245.
\end{prova}

\section{Paraconsistentizing the logic $L_{3}$ of \L ukasiewicz}

The logic $L_{3}$ of \L ukasiewicz initially proposed in \cite{lukas}, and extensively studied in the literature, is characterized by the matrix
$$
L_{3} = \langle \{0, 1/2, 1\}, \{1\}, f_{\neg}, f_{\vee}, f_{\wedge}, f_{\rightarrow} \rangle
$$
such that: 

\medskip
i. $f_\neg(x)= 1- x$;

ii. $f_{\vee}(x,y)= Max \{x,y\}$;

iii. $f_{\wedge}(x,y)= Min\{x,y\}$;

iv. $f_{\rightarrow}(x,y)= Min\{1, (1 - x + y)\}$.

\medskip
These conditions give rise to the following three-valued truth-tables:

$$
\begin{array}{cc|cccc}
x & y & f_{\neg} (y) & f_{\vee} (x,y) & f_{\wedge} (x,y) & f_{\rightarrow} (x,y) \\
\hline
1 & 1 & 0 & 1 & 1 & 1 \\
1 & 1/2 & 1/2 & 1 & 1/2 & 1/2 \\
1 & 0 & 1 & 1 & 0 & 0 \\
1/2 & 1 &  & 1 & 1/2 & 1 \\
1/2 & 1/2 & & 1/2 & 1/2 & 1 \\
1/2 & 0 & & 1/2 & 0 & 1/2 \\
0 & 1 & & 1 & 0 & 1 \\
0 & 1/2 & & 1/2 & 0 & 1 \\
0 & 0 & & 0 & 0 & 1
\end{array}
$$

The consequence relation for $L_{3}$  is defined as follows: $\Gamma \vDash_{L_{3}} \alpha$ iff for all $L_{3}$-valuation $v$ such that $v(\gamma)=1$ for all $\gamma \in \Gamma$, we have $v(\alpha)=1$. (Notice also that $1$ is the only designated value). In this precise way, if we apply a $\mathbb{P}$-transformation on $L_{3}$, we obtain the following
relation:\footnote{In the case of the logic $L_3$, the question of how to paraconsistentize it has been posed by Walter Carnielli (personal communication, 2005) to the second author of this paper.}

\begin{center}
$\Gamma \vDash^{\mathbb{P}}_{L_{3}} \alpha$ iff there exists $\Gamma'\subseteq \Gamma$, $Mod_{L_{3}}(\Gamma') \neq \emptyset$, such that $\Gamma' \vDash_{L_{3}} \alpha$. 
\end{center}

Notice that, by Lemma \ref{consistencia}, we have that $\Gamma$ is $L_{3}$-consistent iff $Mod_{L_{3}}(\Gamma) \neq \emptyset$, that is,
there is a $L_{3}$-valuation $v$ such that $v(\gamma)=1$ for all $\gamma \in \Gamma$. We also have that $\Gamma \vDash_{L_{3}} \alpha$ iff $Mod_{L_{3}}(\Gamma) \subseteq Mod_{L_{3}}
(\alpha)$, by definition.

The paraconsistentization of $L_{3}$ using the standard procedure by means of consistent
sets was first developed in \cite{dias}. Here we implement and improve those results.

\begin{proposicao}
$L_{3}$ satisfies explosion, joint consistency and it has the conjunctive property.
\end{proposicao}

\begin{prova}
(a) If $\{\alpha, \neg \alpha\} \subseteq \Gamma$, we have that $v(\alpha)=1$ iff $v(\neg \alpha)=0$, for all $L_{3}$-valuation $v$. So, $Mod_{L_{3}}(\Gamma) = \emptyset$ and $Cn_{L}(\Gamma) = For$, that is, $\Gamma$ is $L_{3}$-inconsistent and $L_{3}$ is explosive.

(b) In $L_{3}$, for all propositional letters $p$ we have that $v_{1}(p) = 1$ and $v_{2}(\neg p) = 1$ for some $v_{1}, v_{2} \in \mathcal{V}_{L_{3}}$, but $Mod_{L_{3}}(\{p, \neg p\}) = \emptyset$, that is, $\{p, \neg p\}$ is $L_{3}$-inconsistent.

(c) In $L_{3}$, $Cn_{L_{3}}(\{\alpha, \beta\})=Cn_{L_{3}}(\{\alpha \wedge \beta\})$.
\end{prova}

\begin{corolario}
$\mathbb{P}(L_{3})$ is paraconsistent.
\end{corolario}

Let us verify the idempotency of the $\mathbb{P}$-transformation with respect to $L_{3}$. Notice that, in $L_{3}$, we have that $Mod_{L_{3}}(\neg(\alpha \rightarrow \alpha)) = \emptyset$
and, then, $\{\neg(\alpha \rightarrow \alpha)\} \vDash_{L_{3}} \beta$, for all $\beta \in For$.

\begin{corolario}
$\mathbb{P}(L_{3})= \mathbb{P}(\mathbb{P}(L_{3}))$.
\end{corolario}

\begin{prova}
Since $L_{3}$ is a monotonic logic (that is, $A \subset Cn_{L_{3}}(\Gamma)$), if $\Gamma \subseteq For$ is $L_{3}$-consistent, then $Cn_{L_{3}}(\Gamma)=Cn^{\mathbb{P}}_{L_{3}}(\Gamma)$. As $Mod_{L_{3}}(\neg(\alpha \rightarrow \alpha)) = \emptyset$, by the proposition \ref{prop36}, every $\Gamma \subset For$ is $\mathbb{P}(L_{3})$-consistent. Therefore, $Cn^{\mathbb{P}}_{L_{3}}(Cn^{\mathbb{P}}_{L_{3}}(\Gamma))=Cn^{\mathbb{P}}_{L_{3}}(\Gamma)$, for all $\Gamma \subseteq For$.
\end{prova}

\begin{proposicao}
${\mathbb{P}}(L_{3})$ does not satisfy inclusion.
\end{proposicao}

\begin{prova}
It is enough to see that $\{\neg(\alpha \rightarrow \alpha)\} \nvDash^{\mathbb{P}}_{L_{3}} \neg(\alpha \rightarrow \alpha)$.
\end{prova}

\medskip
In $L_{3}$, the \emph{transitivity property} holds: If $\Delta \vDash_{L_{3}} \alpha$ and $\Gamma \vDash_{L_{3}} \delta$, for all $\delta \in \Delta$, then $\Gamma \vDash_{L_{3}} \alpha$.

\begin{proposicao}
\label{transL3}
${\mathbb{P}}(L_{3})$ does not satisfy transitivity.
\end{proposicao}

\begin{prova}
In $L_{3}$, we have that $\{\alpha\} \vDash_{L_{3}} \alpha \vee \beta$. Let $p,q$ be propositional letters. Since $\{p\}$ and $\{\neg p\}$ are $L_{3}$-consistent, it follows that $\{p\} \vDash^{\mathbb{P}}_{L_{3}} p\vee q$ and also $\{\neg p \} \vDash^{\mathbb{P}}_{L_{3}} \neg p$. Thus, if $\Delta = \{p \vee q, \neg p\}$ and $\Gamma=\{p, \neg p\}$, then it follows that $\Gamma \vDash^{\mathbb{P}}_{L_{3}} \delta$, for all $\delta \in \Delta$. On the other hand, if $v \in \mathcal{V}_{L_{3}}$ is such that $v(p \vee q) = v(\neg p)=1$, we have that $v(p) = 0$ and then $v(q)=1$. Therefore, $\Delta \vDash_{L_{3}} q$. Thus, $\Delta \vDash^{\mathbb{P}}_{L_{3}} q$ since that $\Delta$ is $L_{3}$-consistent. Moreover, the only $L_{3}$-consistent subsets of $\Gamma$ are: $\emptyset, \{p\}$ and $\{\neg p\}$. But $\emptyset \nvDash_{L_{3}} q$, $\{p\}\nvDash_{L_{3}} q$ and $\{\neg p\} \nvDash_{L_{3}} q$. So, $\Gamma \nvDash^{\mathbb{P}}_{L_{3}} q$.
\end{prova}

\begin{corolario}
${\mathbb{P}}(L_{3})$ does not satisfy idempotency.
\end{corolario}

\begin{prova}
Consider $\Gamma=\{p,\neg p\}$ ($p$ is a propositional letter). Then, $p\vee q, \neg p \in Cn^{\mathbb{P}}_{L_{3}}(\Gamma)$ and $q \in Cn^{\mathbb{P}}_{L_{3}}(Cn^{\mathbb{P}}_{L_{3}}(\Gamma))$ but $q \notin Cn^{\mathbb{P}}_{L_{3}}(\Gamma)$.
\end{prova}

\begin{proposicao}
${\mathbb{P}}(L_{3})$ is monotonic.
\end{proposicao}

\begin{prova}
In fact, if $L=\langle X, Cn_{L}\rangle$ is a consequence structure, then ${\mathbb{P}}(L)$ is monotonic. See \cite{egs-acl-dd-1}, proposition 3.4, p.245. (Paraconsistentization by consistent sets enforces monotonicity.)
\end{prova}

\medskip
Despite the fact that ${\mathbb{P}}(L_{3})$ does not satisfy transitivity in its full form, a weak form of transitivity holds.

\begin{proposicao}[Weak transitivity]
If $\{\alpha\} \vDash^{{\mathbb{P}}}_{L_{3}} \beta$ and $\{\beta\} \vDash^{{\mathbb{P}}}_{L_{3}} \gamma$, then $\{\alpha\} \vDash^{{\mathbb{P}}}_{L_{3}} \gamma$.
\end{proposicao}

We postponed the proof and consider a previous lemma. We say that $\alpha \in For$ is a $L_{3}$-\emph{contradiction} iff for every $v \in \mathcal{V}_{L_{3}}$ we have that $v(\alpha) = 0$. (We use this definition always when the matrix has a truth-value $0 \in Val$ such that $0 \notin D$.) Moreover, we recall that a formula $\alpha \in For$ is a $L_{3}$-tautology iff for every $v \in \mathcal{V}_{L_{3}}$ we have that $v(\alpha) = 1$. ($1$ is the only designated value.)

\begin{lema}
\label{lemaweak}
It holds that:

i. If $\beta$ is a $L_{3}$-contradiction, then for all $\Gamma \subseteq For$, $\Gamma \nvDash^{{\mathbb{P}}}_{L_{3}} \beta$;

ii. If $\beta$ is a $L_{3}$-tautology and $\{\beta\} \vDash^{{\mathbb{P}}}_{L_{3}} \gamma$, then $\gamma$ is a $L_{3}$-tautology and for every $\Gamma \subseteq For$,  $\Gamma \vDash^{{\mathbb{P}}}_{L_{3}} \gamma$;

iii. If $\{\alpha\} \vDash^{{\mathbb{P}}}_{L_{3}} \beta$, then $\beta$ is a $L_{3}$-tautology or $\{\alpha\}$ is $L_{3}$-consistent and $\{\alpha\} \vDash_{L_{3}} \beta$.
\end{lema}

\begin{prova}
i. Suppose that $\beta$ is a $L_{3}$-contradiction and there is a $\Gamma \subseteq For$ such that $\Gamma \vDash^{{\mathbb{P}}}_{L_{3}} \beta$. So, there is a $\Gamma' \subseteq \Gamma$ with $Mod_{L_{3}}(\Gamma') \neq \emptyset$ such that $\Gamma' \vDash_{L_{3}} \psi$. But, for $v \in Mod_{L_{3}}(\Gamma')$ we
have $v(\beta) = 0$ (contradiction!).

ii. Suppose that $\beta$ is a $L_{3}$-tautology and there is a $v \in \mathcal{V}_{L_{3}}$ such that $v(\gamma) \neq 1$. Since $\{\beta\} \vDash^{{\mathbb{P}}}_{L_{3}} \gamma$, we have $\emptyset \vDash_{L_{3}} \gamma$ or $\{\beta\} \vDash_{L_{3}} \gamma$. In both cases, $\gamma$ has to be a $L_{3}$-tautology. Further,  if $\gamma$ is a $L_{3}$-tautology, then it is clear that $\Gamma \vDash^{{\mathbb{P}}}_{L_{3}} \gamma$, for all $\Gamma \subseteq For$.

iii. Immediate from i. and ii. 
\end{prova}

\medskip
Now, we prove weak transitivity.

\medskip
\begin{provap}
Suppose that $\{\alpha\} \vDash^{{\mathbb{P}}}_{L_{3}} \beta$ and $\{\beta\} \vDash^{{\mathbb{P}}}_{L_{3}} \gamma$. By Lemma (i), $\beta$ is not a $L_{3}$-contradiction. We have two cases:

(a) $\beta$ is a  $L_{3}$-tautology. In this case, since $\{\beta\} \vDash^{{\mathbb{P}}}_{L_{3}} \gamma$, by Lemma (ii), we have that $\{\alpha\} \vDash^{{\mathbb{P}}}_{L_{3}} \gamma$.

(b) $\psi$ is not a  $L_{3}$-tautology. In this case, by Lemma (iii), $\{\alpha\}$ is $L_{3}$-consistent and $\{\alpha\} \vDash_{L_{3}} \beta$. On the other hand, since $\{\beta\} \vDash^{{\mathbb{P}}}_{L_{3}} \gamma$ we have that $\gamma$ is a $L_{3}$-tautology, and we have $\{\alpha\} \vDash^{{\mathbb{P}}}_{L_{3}} \gamma$; or $\{\beta\}$ is $L_{3}$-consistent and  $\{\beta\} \vDash_{L_{3}} \gamma$. Since $\{\alpha\} \vDash_{L_{3}} \beta$ and $\{\beta\} \vDash_{L_{3}} \gamma$, by transitivity in $L_{3}$, we have $Mod_{L_{3}}(\{\alpha\}) \subseteq Mod_{L_{3}}(\{\beta\}) \subseteq Mod_{L_{3}}(\{\gamma\})$ and then $\{\alpha\} \vDash_{L_{3}} \gamma$. Therefore, $\{\alpha\} \vDash^{{\mathbb{P}}}_{L_{3}} \gamma$.
\end{provap}

\medskip
It is a well known fact that the semantic deduction theorem does not hold in $L_{3}$. There is, however, the following version of the deduction theorem:
$$
\Gamma \cup \{\alpha\} \vDash_{L_{3}} \beta \text{ iff } \Gamma \vDash_{L_{3}} \alpha \rightarrow (\alpha \rightarrow \beta). \ \ \ \ \ \ \ (*)
$$
In $\mathbb{P}(L_{3})$ we have only that:

\begin{proposicao}
If $\Gamma \cup \{\alpha\} \vDash^{{\mathbb{P}}}_{L_{3}} \beta$, then $\Gamma \vDash^{{\mathbb{P}}}_{L_{3}} \alpha \rightarrow (\alpha \rightarrow \beta)$.
\end{proposicao}

\begin{prova}
Suppose that $\Gamma \cup \{\alpha\} \vDash^{{\mathbb{P}}}_{L_{3}} \beta$. Then, there is $\Gamma' \subseteq \Gamma \cup \{\alpha\}$,  $L_{3}$-consistent, such that $\Gamma'\vDash_{L_{3}} \beta$. We have
then three cases:

(a) $\alpha \in \Gamma$. In this case, $\Gamma'\subseteq \Gamma$ such that $\Gamma' \cup \{\alpha\}  \vDash_{L_{3}} \beta$ ($L_{3}$ is monotonic). Then, by (*), $\Gamma \vDash_{L_{3}} \alpha \rightarrow (\alpha \rightarrow \beta)$. Since $\Gamma'$ is $L_{3}$-consistent, we have $\Gamma \vDash^{{\mathbb{P}}}_{L_{3}} \alpha \rightarrow (\alpha \rightarrow \beta)$.

(b) $\alpha \notin \Gamma$ and $\alpha \notin \Gamma'$. In this case, we have again that $\Gamma'\subseteq \Gamma$  and the same argument shows that $\Gamma \vDash^{{\mathbb{P}}}_{L_{3}} \alpha \rightarrow (\alpha \rightarrow \beta)$.

(c) $\alpha \notin \Gamma$ and $\alpha \in \Gamma'$. In this case, $\Gamma'-\{\alpha\}\subseteq \Gamma$ is $L_{3}$-consistent and $\Gamma'-\{\alpha\} \cup \{\alpha\} \vDash_{L_{3}} \beta$. Given (*), $\Gamma'-\{\alpha\} \vDash_{L_{3}}  \alpha \rightarrow (\alpha \rightarrow \beta)$. However, $\Gamma'-\{\varphi\} \subseteq \Gamma$ is $L_{3}$-consistent, and we have $\Gamma \vDash^{{\mathbb{P}}}_{L_{3}} \alpha \rightarrow (\alpha \rightarrow \beta)$.
\end{prova}

\medskip
In order to see that the converse of the above proposition is not valid, we consider $\alpha = \beta = \neg(p \rightarrow p)$ ($p$ is a propositional letter). In this case, we have that
$$
\Gamma \vDash^{{\mathbb{P}}}_{L_{3}} \neg(p \rightarrow p) \rightarrow (( \neg(p \rightarrow p) \rightarrow \neg(p \rightarrow p)))
$$
because $\neg(p \rightarrow p) \rightarrow (( \neg(p \rightarrow p) \rightarrow \neg(p \rightarrow p)))$ is a $L_{3}$-tautology. But 
$$
\Gamma \cup \{\neg(p \rightarrow p)\} \nvDash^{{\mathbb{P}}}_{L_{3}} \neg(p \rightarrow p)
$$
because $\neg(p \rightarrow p)$ is a $L_{3}$-contradiction.

\section{Paraconsistentizing the system $G_{3}$ of K. G\"odel}
The logic $G_{3}$ of G\"odel (originally developed in \cite{godel}) is characterized by the matrix
$$
G_{3} = \langle \{0, 1/2, 1\}, \{1\}, f_{\neg}, f_{\vee}, f_{\wedge}, f_{\rightarrow} \rangle
$$
such that: 

\medskip
i. $f_\neg(x)= \left\{
\begin{array}{ll}
1 & \text{if } x = 0 \\
0 & \text{if } x \neq 0;
\end{array}
\right.
$

ii. $f_{\vee}(x,y)= Max \{x,y\}$;

iii. $f_{\wedge}(x,y)= Min\{x,y\}$;

iv. $f_{\rightarrow}(x,y)= \left\{
\begin{array}{ll}
1 & \text{if } x \leq y \\
y & \text{if } x > y;
\end{array}
\right.
$.

\medskip
These conditions give rise to the following three-valued truth-tables:
$$
\begin{array}{cc|cccc}
x & y & f_{\neg} (y) & f_{\vee} (x,y) & f_{\wedge} (x,y) & f_{\rightarrow} (x,y) \\
\hline
1 & 1 & 0 & 1 & 1 & 1 \\
1 & 1/2 & 0 & 1 & 1/2 & 1/2 \\
1 & 0 & 1 & 1 & 0 & 0 \\
1/2 & 1 &  & 1 & 1/2 & 1 \\
1/2 & 1/2 & & 1/2 & 1/2 & 1 \\
1/2 & 0 & & 1/2 & 0 & 0 \\
0 & 1 & & 1 & 0 & 1 \\
0 & 1/2 & & 1/2 & 0 & 1 \\
0 & 0 & & 0 & 0 & 1
\end{array}
$$

The difference between $G_{3}$ and $L_{3}$ remains in the definition of $f_{\rightarrow}$ and $f_{\neg}$. In $G_{3}$, $f_{\rightarrow}(1/2,0) = 0$, while in $L_{3}$, $f_{\rightarrow}(1/2,0) = 1/2$. On the other hand, $f_{\neg}(1/2) = 0$ in $G_{3}$, and $f_{\neg}(1/2) = 1/2$ in $L_{3}$. The consequence relations $\vDash_{G_{3}}$ and $\vDash^{\mathbb{P}}_{G_{3}}$ are:
$$
\Gamma \vDash_{G_{3}} \alpha \Leftrightarrow Mod_{G_{3}}(\Gamma) \subseteq Mod_{G_{3}}(\alpha),
$$
and, by Lemma \ref{consistencia},
$$
\Gamma \vDash^{\mathbb{P}}_{G_{3}} \alpha \Leftrightarrow \text{ there exists } \Gamma' \subseteq \Gamma, Mod_{G_{3}}(\Gamma') \neq \emptyset, \text{ such that } \Gamma' \vDash_{G_{3}} \alpha.
$$

It is easy to see that $G_{3}$ satisfies explosive property, because if $\{\alpha, \neg \alpha\} \subseteq Cn_{G{3}}(\Gamma)$, then $Mod_{G_{3}}(\Gamma) \subseteq Mod_{G_{3}}(\{\alpha,\neg\alpha\}) = \emptyset$ and, therefore, $\Gamma \vDash_{G_{3}} \beta$ for all $\beta \in For$. Moreover, $G_{3}$ satisfies joint consistency and conjunctive property. Thus, we have:

\begin{proposicao}
$\mathbb{P}(G_{3})$ is paraconsistent.
\end{proposicao}

On the other hand, we have that $\alpha \wedge \neg \alpha$ is a $G_{3}$-contradiction and then $\{\alpha \wedge \neg \alpha\}$ is $G_{3}$-inconsistent. So, by Proposition \ref{prop36}, we have:

\begin{proposicao}
$\mathbb{P}(\mathbb{P}(G_{3})) = \mathbb{P}(G_{3})$.
\end{proposicao}

\begin{proposicao}
$\mathbb{P}(G_{3})$ does not satisfy inclusion.
\end{proposicao}

\begin{prova}
$\{\alpha \wedge \neg \alpha\} \nvDash^{\mathbb{P}}_{G_{3}} \alpha \wedge \neg \alpha$.
\end{prova}

\begin{proposicao}
$\mathbb{P}(G_{3})$ does not satisfy transitivity.
\end{proposicao}

\begin{prova}
Analogous to the proof of proposition 8.
\end{prova}

\begin{corolario}
${\mathbb{P}}(G_{3})$ does not satisfy idempotency.
\end{corolario}

Let us consider the weak transitivity. In this case, the Lemma \ref{lemaweak} remains valid for ${\mathbb{P}}(G_{3})$. Thus, we have:

\begin{proposicao}
In ${\mathbb{P}}(G_{3})$, the weak form of transitivity is valid.
\end{proposicao}

\begin{prova}
Analogous to the proof of proposition 11.
\end{prova}

\medskip
It remains to consider the deduction theorem. It is well known that the $G_{3}$-matrix ``involves all assertions of intuitionistic propositional calculus. G\"odel's axiomatization coincides with the usual axiom system of the intuitionistic calculus with
$$
(\neg \alpha \rightarrow \beta) \rightarrow (((\beta \rightarrow \alpha) \rightarrow \beta) \rightarrow \beta))
$$
added'' (see Bolc \& Borowik \cite{bolc}, p 84). Therefore, the deduction theorem is valid in $G_{3}$.

In $\mathbb{P}(G_{3})$ we have only that:

\begin{proposicao}
If $\Gamma \cup \{\alpha\} \vDash^{{\mathbb{P}}}_{G_{3}} \beta$, then $\Gamma \vDash^{{\mathbb{P}}}_{G_{3}} \alpha \rightarrow \beta$.
\end{proposicao}

\begin{prova}
Suppose that $\Gamma \cup \{\alpha\} \vDash^{{\mathbb{P}}}_{G_{3}} \beta$. Then, there is $\Gamma' \subseteq \Gamma \cup \{\alpha\}$,  $G_{3}$-consistent, such that $\Gamma'\vDash_{G_{3}} \beta$. We have
then three cases:

(a) $\alpha \in \Gamma$;

(b) $\alpha \notin \Gamma$ and $\alpha \notin \Gamma'$.

In these cases, $\Gamma'\subseteq \Gamma$ and $\Gamma'\vDash_{G_{3}} \beta$. By monotonicity, $\Gamma' \cup \{\alpha\}\vDash_{G_{3}} \beta$ and, since deduction theorem is valid in $G_{3}$, we have $\Gamma'\vDash_{G_{3}} \alpha \rightarrow \beta$. Thus, $\Gamma \vDash^{{\mathbb{P}}}_{G_{3}} \alpha \rightarrow \beta$.

(c) $\alpha \notin \Gamma$ and $\alpha \in \Gamma'$. In this case, let $\Gamma^{*} = \Gamma' - \{\alpha\}$. Since $\Gamma^{*} \subseteq \Gamma'$, $\Gamma^{*}$ is $G_{3}$-consistent and $\Gamma' = \Gamma^{*} \cup \{\alpha\} \vDash_{G_{3}} \beta$. By deduction theorem for $G_{3}$, we have $\Gamma^{*} \vDash_{G_{3}} \alpha \rightarrow\beta$. Therefore, $\Gamma \vDash^{{\mathbb{P}}}_{G_{3}} \alpha \rightarrow \beta$.
\end{prova}

\medskip
The converse of proposition above is not valid. It is enough to see that $(\alpha \wedge \neg \alpha) \rightarrow (\alpha \wedge \neg \alpha)$ is a $G_{3}$-tautology and we have $\vDash^{{\mathbb{P}}}_{G_{3}} (\alpha \wedge \neg \alpha) \rightarrow (\alpha \wedge \neg \alpha)$, but $\{\alpha \wedge \neg \alpha \} \nvDash^{{\mathbb{P}}}_{G_{3}} \alpha \wedge \neg \alpha$.

\section{Paraconsistentizing the system $K_{3}$ of S. Kleene}
The system $K_{3}$ of Kleene (see \cite{kleene}) is very similar to $L_{3}$. There is just one difference in the definition of $f_{\rightarrow}$. In $K_{3}$, $f_{\rightarrow}(1/2,1/2) = 1/2$. Thus, $K_{3}$ is characterized by the matrix
$$
K_{3} = \langle \{0, 1/2, 1\}, \{1\}, f_{\neg}, f_{\vee}, f_{\wedge}, f_{\rightarrow} \rangle
$$
such that $f_{\neg}, f_{\vee}, f_{\wedge}$ and $f_{\rightarrow}$ is given by the following truth tables:
$$
\begin{array}{cc|cccc}
x & y & f_{\neg} (y) & f_{\vee} (x,y) & f_{\wedge} (x,y) & f_{\rightarrow} (x,y) \\
\hline
1 & 1 & 0 & 1 & 1 & 1 \\
1 & 1/2 & 1/2 & 1 & 1/2 & 1/2 \\
1 & 0 & 1 & 1 & 0 & 0 \\
1/2 & 1 &  & 1 & 1/2 & 1 \\
1/2 & 1/2 & & 1/2 & 1/2 & 1/2 \\
1/2 & 0 & & 1/2 & 0 & 1/2 \\
0 & 1 & & 1 & 0 & 1 \\
0 & 1/2 & & 1/2 & 0 & 1 \\
0 & 0 & & 0 & 0 & 1
\end{array}
$$

Let $\alpha(p_{1},...p_{n})$ be a formula such that $\alpha$ has $p_{1},...p_{n}$ as propositional letters. Let $v \in \mathcal{V}_{K_{3}}$ be a $K_{3}$-valuation such that $V(p_{1}) = ... = v(p_{n}) = 1/2$. Then, it is easy to see that $v(\alpha) = 1/2$. We conclude that: \emph{In $K_{3}$, we have neither tautologies nor contradictions}.

We recall that $1$ is the unique designated value and if $\Gamma \subseteq For$, then
$$
Mod_{K_{3}}(\Gamma) = \{v \in \mathcal{V}_{K_{3}}: v(\gamma) = 1 \text{ for all } \gamma \in \Gamma\}.
$$
The consequence relations are defined as usual: 
$$
\Gamma \vDash_{K_{3}} \alpha \Leftrightarrow Mod_{K_{3}}(\Gamma) \subseteq Mod_{K_{3}}(\alpha),
$$
and
$$
\Gamma \vDash^{\mathbb{P}}_{K_{3}} \alpha \Leftrightarrow \text{ there exists } \Gamma' \subseteq \Gamma, K_{3}\text{-consistent, such that } \Gamma' \vDash_{K_{3}} \alpha.
$$
But, by Lemma \ref{consistencia}, in $K_{3}$ we have: $\Gamma \subseteq For$ is $K_{3}$-consistent iff $Mod_{K_{3}}(\Gamma) \neq \emptyset$.

Although $K_{3}$ has no contradictions, the set $\{\alpha \wedge \neg\alpha\}$ is $K_{3}$-inconsistent (there is no $K_{3}$-valuation $v$ such that $v(\alpha \wedge \neg\alpha) = 1$), that is, $Mod_{K_{3}}(\{\alpha \wedge \neg\alpha\}) = \emptyset$ and $Cn_{K_{3}}(\{\alpha \wedge \neg\alpha\}) = For$.

Since $K_{3}$ was defined by means of a matrix, $K_{3}$ is a normal logic, that is, $K_{3}$ satisfies inclusion, monotonicity and idempotency. Moreover, $K_{3}$ satisfies explosive property (if $\{\alpha,\neg\alpha\} \subseteq Cn_{K_{3}}(\Gamma)$, then $Mod_{K_{3}}(\Gamma) = \emptyset$ and $Cn_{K_{3}}(\{\Gamma) = For$), joint consistency (by propositional letters) and conjunctive property ($Cn_{K_{3}}(\{\alpha, \beta\}) = Cn_{K_{3}}(\{\alpha \wedge \beta\})$).

\begin{proposicao}
$\mathbb{P}(K_{3})$ is paraconsistent.
\end{proposicao}

Since $\{\alpha \wedge \neg\alpha\}$ is $K_{3}$-inconsistent, by Proposition \ref{prop36}, $\mathbb{P}(K_{3})$ has no $\mathbb{P}(K_{3})$-inconsistent sets, and we have:

\begin{proposicao}
$\mathbb{P}(\mathbb{P}(K_{3})) = \mathbb{P}(K_{3})$.
\end{proposicao}

Since $K_{3}$ has no tautology, $Mod_{K_{3}}(\emptyset) = \mathcal{V}_{K_{3}}$ and $Cn_{K_{3}}(\emptyset) = \emptyset$. Thus, the unique $K_{3}$-consistent subset of $\{\alpha \wedge \neg\alpha\}$ is $\emptyset$. So, $\{\alpha \wedge \neg\alpha\} \nvDash^{\mathbb{P}}_{K_{3}} \alpha \wedge \neg\alpha$, and we have:

\begin{proposicao}
$\mathbb{P}(K_{3})$ does not satisfy inclusion.
\end{proposicao}

The same counter-example that shows that $\mathbb{P}(L_{3})$ does not satisfy transitivity (Proposition \ref{transL3}) can be used in $\mathbb{P}(K_{3})$. For propositional letters $p$ and $q$, $\{p \vee q, \neg p\} \vDash^{\mathbb{P}}_{K_{3}} q$, $\{p,\neg p\} \vDash^{\mathbb{P}}_{K_{3}} p \vee q$, $\{p,\neg p\} \vDash^{\mathbb{P}}_{K_{3}} p \vee \neg p$, but $\{p,\neg p\} \nvDash^{\mathbb{P}}_{K_{3}} q$. Thus, we have:

\begin{proposicao}
$\mathbb{P}(K_{3})$ does not satisfy transitivity.
\end{proposicao}

\begin{corolario}
${\mathbb{P}}(K_{3})$ does not satisfy idempotency.
\end{corolario}

Now, we consider the weak transitivity.

\begin{lema}
In ${\mathbb{P}}(K_{3})$, if $\{\alpha\} \vDash^{\mathbb{P}}_{K_{3}} \beta$, then $\{\alpha\}$ is $K_{3}$-consistent and $\{\alpha\} \vDash_{K_{3}} \beta$.
\end{lema}

\begin{prova}
Suppose that $\{\alpha\} \vDash^{\mathbb{P}}_{K_{3}} \beta$. Then, there exists $\Gamma \subseteq \{\alpha\}$, $K_{3}$-consistent, such that $\Gamma \vDash_{K_{3}} \beta$. But $\Gamma$ cannot be $\emptyset$ because $Mod_{K_{3}}(\emptyset) = \mathcal{V}_{K_{3}}$ and, in this case, $Mod_{K_{3}}(\beta) = \mathcal{V}_{K_{3}}$ and $\beta$ would be a tautology. Therefore, $\{\alpha\}$ has to be $K_{3}$-consistent and $\{\alpha\} \vDash_{K_{3}} \beta$.
\end{prova}

\begin{proposicao}
$\mathbb{P}(K_{3})$ satisfies weak transitivity.
\end{proposicao}

\begin{prova}
Suppose that $\{\alpha\} \vDash^{\mathbb{P}}_{K_{3}} \beta$ and $\{\beta\} \vDash^{\mathbb{P}}_{K_{3}} \gamma$. By Lemma, we have that $\{\alpha\}$ and $\{\beta\}$ are $K_{3}$-consistent, $\{\alpha\} \vDash_{K_{3}} \beta$ and $\{\beta\} \vDash_{K_{3}} \gamma$. By transitivity in $K_{3}$, we have $\{\alpha\} \vDash_{K_{3}} \gamma$. Therefore, $\{\alpha\} \vDash^{\mathbb{P}}_{K_{3}} \gamma$.
\end{prova}

\medskip
It is easy to see that deduction theorem is not valid in $K_{3}$. In this system, we have $\{\alpha\} \vDash_{K_{3}} \alpha$ (by inclusion), but $\nvDash_{K_{3}} \alpha \rightarrow \alpha$. The same example, for propositional letters, shows that deduction theorem is not valid in $\mathbb{P}(K_{3})$ also. Below we present a table containing a logic and its paraconsistentized version.
The reader can compare properties which the initial logic has and what is preserved (or lost)
in the paraconsistentized version of it. 

\begin{center}
\textbf{Summary of results}
\end{center}

$$
\begin{array}{l|cccccc}
 & L_{3} & \mathbb{P}(L_{3}) & G_{3} & \mathbb{P}(G_{3}) & K_{3} & \mathbb{P}(K_{3}) \\
\hline
\text{explosive property} & \checkmark & \times & \checkmark & \times & \checkmark & \times \\
\text{joint consistency} & \checkmark & \checkmark & \checkmark & \checkmark & \checkmark & \checkmark \\
\text{conjunctive property} & \checkmark & \times & \checkmark & \times & \checkmark & \times \\
\text{paraconsistent} & \times & \checkmark & \times & \checkmark & \times & \checkmark \\
\text{inconsistent sets} & \checkmark & \times & \checkmark & \times & \checkmark & \times \\
\mathbb{P}(\mathbb{P}(L)) = \mathbb{P}(L) & \checkmark & \checkmark & \checkmark & \checkmark & \checkmark & \checkmark \\
\text{inclusion} & \checkmark & \times & \checkmark & \times & \checkmark & \times \\
\text{monotonicity} & \checkmark & \checkmark & \checkmark & \checkmark & \checkmark & \checkmark \\
\text{idempotency} & \checkmark & \times & \checkmark & \times & \checkmark & \times \\
\text{transitivity} & \checkmark & \times & \checkmark & \times & \checkmark & \times \\
\text{weak transitivity} & \checkmark & \checkmark & \checkmark & \checkmark & \checkmark & \checkmark \\
\text{modus ponens} & \checkmark & \times & \checkmark & \times & \checkmark & \times \\
\text{full deduction theorem} & \times & \times & \checkmark & \times & \times & \times \\
\text{modified full deduction theorem} & \checkmark & \checkmark & \checkmark & \times & \times & \times \\
\text{weak deduction theorem } (\Rightarrow) & \times & \times & \checkmark & \checkmark & \times & \times \\
\text{modified weak deduction theorem } (\Rightarrow) & \checkmark & \checkmark & \checkmark & \checkmark & \times & \times \\
\end{array}
$$
A tick ($\checkmark$) means \emph{yes}; a cross ($\times$) means \emph{no}.

\section{The logic of paradox}

Now, we investigate the interesting case of what happens when
the $\mathbb{P}$-transformation is applied to an already paraconsistent
logic. Although it seems unuseful to apply it to a logic which is already
paraconsistent, it is interesting to check whether the paraconsistent logic 
remains the same or not. The case of the
logic of paradox (designed by Priest in \cite{priest0}) seems to be of special
interest as it is paraconsistent and many-valued. The logic of
paradox has very special and unique characteristics which makes it a rather
complicated case. So, in this section, we study a paraconsistentization of 
a logic that is already paraconsistent. Let's check it how to paraconsistentize
the logic of paradox $LP$.

The logic of paradox ($LP$) of G. Priest is characterized by the logical matrix
$$
LP = \langle \{t,p,f\}, \{t,p\}, f_{\neg}, f_{\vee}, f_{\wedge}, f_{\rightarrow} \rangle
$$
such that $\{t,p\}$ is the set of designated truth-values and $f_{\neg}, f_{\vee}, f_{\wedge}$, $f_{\rightarrow}$ are functions in $\{t,p,f\}$ given by the following truth-tables
$$
\begin{array}{cc|cccc}
x & y & f_{\neg} (y) & f_{\vee} (x,y) & f_{\wedge} (x,y) & f_{\rightarrow} (x,y) \\
\hline
t & t & f & t & t & t \\
t & p & p & t & p & p \\
t & f & t & t & f & f \\
p & t &  & t & p & t \\
p & p & & p & p & p \\
p & f & & p & f & p \\
f & t & & t & f & t \\
f & p & & p & f & t \\
f & f & & f & f & t
\end{array}
$$

The only difference between $LP$ and the system $K_{3}$ of S. Kleene is the set of designated values.
We introduce two specific $LP$-valuations that will be useful to prove some results.
So, let $v^{p}$ be the $LP$-valuation such that $v^{p}(q) = p$, for all $q \in Prop$. In this case, using the truth-tables, we have that $v^{p}(\alpha) = p$, for all $\alpha \in For$. Therefore, since $p$ is a designated value, it holds that:

\begin{proposicao}
\label{p-model}
In $LP$, $v^{p}$ is a $LP$-model of all subsets $\Gamma \subseteq For$, that is, all sets of $LP$-formulas have a $LP$-model.
\end{proposicao}

On the other hand, let $v^{t}$ be the $LP$-valuation such that $v^{t}(q) = t$, for all $q \in Prop$. Again, using the truth-tables, we have that $v^{t}(\alpha) = t \text{ or } f$, for all $\alpha \in For$. Thus, we have:

\begin{proposicao}
\label{no-p}
In $LP$, there is no formula $\alpha$ such that $v(\alpha) = p$ for all $v \in \mathcal{V}_{LP}$.
\end{proposicao}

Therefore, we have two distinguished (canonical) $LP$-valuations, $v^{P}$ and $v^{t}$, that will imply some consequences in the logic of the system $LP$. 

Let $\alpha \in For$. We define: $\alpha$ is a $LP$-\emph{contradiction} iff $v(\alpha) = f$ for all $v \in \mathcal{V}_{LP}$, that is, $Mod_{LP}(\alpha) = \emptyset$.

\medskip
It is easy to see that:

\begin{proposicao}
In $LP$, there is no $LP$-contradiction.
\end{proposicao}

\begin{prova}
By Proposition \ref{p-model}, for all $\alpha \in For$, $v^{p}(\alpha) = p$.
\end{prova}

\begin{proposicao}
If $\alpha$ is a $LP$-tautology, then there exists $v \in \mathcal{V}_{LP}$ such that $v(\alpha) = t$.
\end{proposicao}

\begin{prova}
Consequence of Proposition \ref{no-p}.
\end{prova}

\medskip
In $LP$, the set $For$, by inclusion, is $LP$-inconsistent. But there are many other $LP$-inconsistent sets.

Let $Var(\Gamma)$ be the set of all propositional variables occurring in $\Gamma$. If $Var(\Gamma) \varsubsetneq Prop$, the for all $q \notin Var(\Gamma)$ we have that $\Gamma \nvDash_{LP} q$. Because, if $v \in \mathcal{V}_{LP}$ is such that $v(r) = p$ for all $r \in Var(\Gamma)$ and $v(q) = f$, then we have $v \in Mod_{LP}(\Gamma)$ and $v \notin Mod_{LP}(q)$, that is, $\Gamma \nvDash_{LP} q$.

Thus, we have the following results:

\begin{proposicao}
If $Var(\Gamma) \varsubsetneq Prop$, then $\Gamma$ is $LP$-consistent.
\end{proposicao}

\begin{corolario}
\label{fin-con}
If $\Gamma$ is finite, then $\Gamma$ is $LP$-consistent.
\end{corolario}

On the other hand, we can produce an infinite $LP$-inconsistent sets. 

\begin{proposicao}
If $\Gamma \subseteq Prop$, then $For - \Gamma$ is $LP$-inconsistent.
\end{proposicao}

\begin{prova}
Since, in $LP$, we have $\{\alpha \wedge \beta\} \vDash_{LP} \alpha$, for all $q \in \Gamma$ we have $(q \wedge q) \in For - \Gamma$ and $For - \Gamma \vDash_{LP} q$. Thus, $For - \Gamma$ is $LP$-inconsistent.
\end{prova}

\medskip
Since, in $LP$, every subset $\Gamma \subseteq For$ has a $LP$-model (Proposition \ref{p-model}), it does not hold the property: $\Gamma$ is $LP$-consistent iff $Mod_{LP}(\Gamma) \neq \emptyset$. Observe that $v^{p} \in Mod_{LP}(For)$, but $For$ is $LP$-inconsistent. So, we have a problem to characterized the $LP$-consistent sets in terms of $LP$-valuations.

\begin{teorema}
$\Gamma$ is $LP$-inconsistent iff $v^{p}$ is the unique $LP$-model of $\Gamma$, that is, $Mod_{LP}(\Gamma) = \{v^{p}\}$.
\end{teorema}

\begin{prova}
($\Rightarrow$) Let $\Gamma$ be a $LP$-inconsistent set and suppose that there is a $LP$-model $v$ such that, for some propositional variable $q$, we have that $v(q) = t \text{ or } f$. If $v(q) = f$, then $\Gamma \nvDash_{LP} q$ and $\Gamma$ is $LP$-consistent (contradiction!). If $v(q) = q$, then $v(\neg q) = f$, $\Gamma \nvDash \neg q$ and, again, $\Gamma$ is $LP$-consistent (contradiction!). So, $Mod_{LP}(\Gamma) = \{v^{p}\}$.

($\Leftarrow$) Suppose that $v^{p}$ is the unique $LP$-model of $\Gamma$. Since $v^{p}(\alpha) = p$, for all $\alpha \in For$ (Proposition \ref{p-model}), we have that $\Gamma \vDash_{LP} \alpha$ for all $\alpha \in For$, that is, $\Gamma$ is $LP$-inconsistent.
\end{prova}

\begin{proposicao}
$LP$ is not compact.
\end{proposicao}

\begin{prova}
Consequence of Corolary \ref{fin-con}.
\end{prova}

\medskip
On the other hand, it is well-known that there is an adequate (sound and complete) axiomatic system $S$ for $LP$ (cf. \cite{priest0}). In this way, we have that: $\Gamma \vDash_{LP} \alpha$ iff there is a finite $\Gamma' \subseteq \Gamma$ such that $\Gamma \vDash_{LP} \alpha$.

\begin{teorema}
$LP$ is a fixed-point of $\mathbb{P}$-transformation. That is, $\mathbb{P}(LP) = LP$.
\end{teorema}

\begin{prova}
Suppose that $\Gamma \vDash_{LP} \alpha$. By the observation above we have that there is a finite $\Gamma' \subseteq \Gamma$ such that $\Gamma \vDash_{LP} \alpha$. By corollary \ref{fin-con}, $\Gamma'$ is $LP$-consistent. Therefore, $\Gamma \vDash^{\mathbb{P}}_{LP} \alpha$ and we have $\mathbb{P}(LP) = LP$.
\end{prova}

\section{Conclusion}

The procedure used here to paraconsistentize some many-valued logics can be
applied to a wide range of logics independent of the fact that they
are many-valued or not. Other procedures of paraconsistentization could
also be developed generating different systems (cf. \cite{egs-acl-dd-3}). 
We have considered only three-valued logics for uniformity, assuming
that  $K_{3}$ and $LP$ are mainly three-valued logics. In the cases
of  $L_{3}$ and  $G_{3}$, the route towards generalization by means of $n$-valued ($n>3$)
systems is rather straightforward, and we have not considered it, though it is virtually
possible to repeat the same steps to get the whole hierarchies. 

G\"odel in his note \cite{godel} has shown that there is an infinite hierarchy of systems
between Heyting's intuitionistic logic $H$ and classical logic, as we mentioned
in the introduction. This infinite hierarchy gives rise to G\"odel intermediate logics. We have produced a paraconsistentization of $G_3$. But the whole hierarchy of G\"odel's
logic could be paraconsistentized and a natural question would be: are these
paraconsistentized intermediate logics between a paraconsistentized version of
$H$ and a paraconsistentized version of classical logic? This topic could be
developed in connection with a paraconsistentization of Heyting's intuitionistic logic.

In \cite{egs-acl-dd-1}, paraconsistentization has been studied in the realm
of abstract logic. In \cite{egs-acl-dd-2}, it has been applied to axiomatic
formal systems and valuation structures, which are still general and abstract.
In this paper, specific many-valued systems were studied applying techniques and methodologies
of paraconsistentization. Departing from four major systems of three-valued logics,
we have explored properties preserved or lost when a paraconsistent transformation
is applied. Due to the broad domain concerning applications of many-valued logics, we can say  that wherever these systems are used to model some concrete and specific situation
and in which types and forms of contradictions appear, our transformed paraconsistent
many-valued systems can play a role, beyond the pure theoretical grounds established
in this paper. In particular, we have suggested that the paraconsistentized many-valued systems explored here can be used and applied to model partial justifications we find in sciences. This occurs especially in virtue of the fact that given a sentence $p$ encoded with its natural weak justification and given also the negation
of $p$ (also encoded with its regular partial justification), these two formulas can be handled inside of the paraconsistentized many-valued systems without explosion. So, all these facts enlarge our understanding and comprehension on the nature, structure, relevance and applicability 
of many-valued logics.


\newpage
\begin{thebibliography}{99}




\bibitem{amgoud}
Amgoud, L.; Besnard, P. (2010). A Formal Analysis of Logic-Based
Argumentation Systems. In: \emph{Scalable Uncertainty Management}, A. Deshpande and A. Hunter (editors). Springer-Verlag Berlin Heidelberg, pp.42-55.

\bibitem{arieli2}
Arieli, O.; Borg, A.; Stra\ss er, C. (2017). Argumentative approaches to reasoning with consistent subsets of premises. In: \emph{Proceedings of the 30th  International  Conference  on  Industrial,  Engineering,  Other  Applications  of  Applied  Intelligent  Systems}, S.Benferhat  \emph{et  alia} (editors), Springer, pp. 455-465.

\bibitem{arieli-avron}
Arieli, O; Avron, A. (2015).  Three-Valued paraconsistent propositional logics.
In: \emph{New Directions in Paraconsistent Logic}, edited by J-Y B\'eziau,
M. Chakraborty and S. Dutta, pp. 91-129. New Delhi: Springer India. 

\bibitem{arieli-avron-zamanski}
Arieli, O; Avron, A, Zamanski, A. (2015). What Is an Ideal Logic for Reasoning with Inconsistency?, In: \emph{Proceedings of the Twenty-Second International Joint Conference on Artificial Intelligence}, pp. 706-711.

\bibitem{avron}
Avron, A.; Lev, I.(2001). A Formula-Preferential Base for
Paraconsistent and Plausible Reasoning Systems. In: \emph{Proceedings of the Workshop on Inconsistency in Data and Knowledge}, pp.60-70. 

\bibitem{benferhat}
Benferhat, S.; Dubois, D.; Prade, H.(1995). How to infer from
inconsistent beliefs without revising? In: \emph{Proceedings of the 14th international joint conference on Artificial intelligence}, vol. 2, pp. 1449-1455   

\bibitem{bensusan-greg}
Bensusan, H; Carneiro, G. (2020). Paraconsistentization through antimonotonicity: towards a logic of supplement.
In: \emph{Abstract Consequence and Logics: Essays in Honor of Edelcio G. de Souza}, edited by A. Costa-Leite,
London: College Publications, pp.263-274.

\bibitem{bolc}
Bolc, L., Borowik P. (1992). \emph{Many Valued Logics 1: theoretical foundations}. Berlin: Springer-Verlag.

\bibitem{brown}
Brown, B.; Priest, G. (2004). Chunk and permeate, a paraconsistent
inference strategy. Part I: The infinitesimal calculus. \emph{Journal of
Philosophical Logic}, 33(4), pp. 379-388.

\bibitem{coniglio1}
Coniglio, M; Esteva, F; Gispert, J; Godo, L. (2019). Maximality in finite-valued \L ukasiewicz logics defined by order filters. \emph{Journal of Logic and Computation},
29(1), pp. 125-156.

\bibitem{costaleite}
Costa-Leite, A. (2007). \emph{Interactions of metaphysical and epistemic concepts}. PhD Thesis, University of
Neuch\^atel, Switzerland.

\bibitem{costaleite3}
Costa-Leite, A. (2014). Lógicas da justificação e quase-verdade. \emph{Principia: an international journal of epistemology}, 18 (2), p.175-186;

\bibitem{costaleite4}
Costa-Leite, A. (2018). O problema das justificações parciais. \emph{Revista de Filosofia Moderna e Contemporânea}, 6(2), pp. 95-104;


\bibitem{dacosta}
da Costa, N. (1974). On the theory of inconsistent formal systems. \emph{Notre Dame Journal
of Formal Logic}, 15, pp.497-510.

\bibitem{dacosta11}
da Costa, N; Alves, E.H. (1977). A semantical analysis of the calculi $C_{n}$. \emph{Notre Dame Journal of Fomal Logic}, 18, pp.621-630.

\bibitem{dacosta-vernengo}
da Costa, N; Vernengo, R. J. (1996). Sobre algunas lógicas paraclásicas y el análisis
del razonamiento jurídico. \emph{Doxa: Cuadernos de Filosofia del Derecho}, 19, pp.183-200.


\bibitem{dacosta1}
da Costa, N. (1999). \emph{O conhecimento científico}, 2a. edição, São Paulo: Discurso Editorial. 


\bibitem{egs-acl-dd-1}
de Souza, E. G; Costa-Leite, A; Dias, D.H.B. (2016). On a paraconsistentization functor in the category of 
consequence structures. \emph{Journal of Applied Non-Classical Logics}, 26(3), pp. 240-250.

\bibitem{egs-acl-dd-2}
de Souza, E. G; Costa-Leite, A; Dias, D.H.B. (2019). Paradeduction in Axiomatic Formal Systems. \emph{Logique
et Analyse}, v. 246, pp. 161-176.

\bibitem{egs-acl-dd-3}
de Souza, E. G; Costa-Leite, A; Dias, D.H.B. (2021). Paraconsistent orbits of logics. \emph{Logica
Universalis}, 15(3), pp. 271-289.

\bibitem{dias}
Dias, D. H. B. (2019). Paraconsistentization of Logics. PhD thesis (In Portuguese),
Universidade de S\~ao Paulo. 

\bibitem{desouza5}
de Souza, E. G. (1998). Remarks on paraclassical logic. \emph{Boletim da Sociedade
Paranaense de Matemática}, 18, pp.107-112.

\bibitem{epstein}
Epstein,  R.  (2001).  Propositional  logics:   the  semantic
foundations  of logic. Wadsworth/Thomson Learning.

\bibitem{godel}
G\"odel, K.(2001). An interpretation of the intuitionistic propositional calculus. \emph{Collected
works}, vol.1, publications 1929-1936. New York: Oxford University Press.

\bibitem{gottwald}
Gottwald, S. (2001). \emph{A Treatise of Many-Valued Logics}. London: Research Studies Press. 

\bibitem{grant-sub}
Grant, J; Subrahmanian,  V.S. The Optimistic and Cautious Semantics for
Inconsistent Knowledge Bases. \emph{Acta Cybernetics}, 12(1), 1995.

\bibitem{jaskowski}
Ja\'skowski, S. (1999). Propositional calculus for contradictory deductive systems. 
\emph{Logic and Logical Philosophy}, 7, pp. 35-56.

\bibitem{kleene}
Kleene, S. C. (2009). \emph{Introduction to Metamathematics}. New York: Ishi Press
International.  

\bibitem{lukas}
Lukasiewicz, J. (1970). On three-valued logic. In: \emph{Jan \L ukasiewicz: selected
works}, edited by L. Borkowski. Amsterdam: North-Holland Publishing. 

\bibitem{malinowski}
Malinowski, G. (1993). \emph{Many-Valued Logics}. New York: Oxford University Press.


\bibitem{payette}
Payette, G. (2009). Preserving logical structure. \emph{On preserving: essays on
preservationism and paraconsistent logic}, pp.105-143, Toronto: University Toronto Press.

\bibitem{priest0}
Priest, G. (1979). The logic of paradox. \emph{Journal of Philosophical Logic} 8, pp.219-241.

\bibitem{priest1}
Priest, G. (2008). \emph{An Introduction to Non-Classical Logic}. Cambridge: Cambridge University Press.  

\bibitem{rescher-manor}
Rescher, N., Manor, R. (1970). On inferences from inconsistent premisses. \emph{Theory and decision}, 1(2), pp.179-217.

\bibitem{subrahmanian}
Subrahmanian, V.; Amgoud, L. (2007). A General Framework for
Reasoning about Inconsistency. In: \emph{Proceedings of the Twentieth International Joint Conference on Artificial Intelligence}, pp.599-604.

\bibitem{tarski}
Tarski, A. (1930). On some fundamental concepts of metamathematics. In:\emph{ Logic, Semantic, Metamathematics}. Second Edition. John Corcoran (ed.). Hackett Publishing Company. 1983. 

 \bibitem{tarski1}
Tarski, A. (1930). Fundamental concepts of the methodology of deductive sciences. In:\emph{Logic, Semantic, Metamathematics}. Second Edition. John Corcoran (ed.). Hackett Publishing Company. 1983. 

\bibitem{vandeputte}
van de Putte, F. (2013). Default assumptions and selection functions: a
generic framework for non-monotonic logics. In: \emph{MICAI 2013: Advances in Artificial Intelligence and Its Applications}, pp. 54-67.


\end{thebibliography}
\end{document}